\documentclass[a4paper,11pt]{article}

\usepackage[plainpages=false]{hyperref}
\usepackage{amsfonts,latexsym,rawfonts,amsmath,amssymb,amsthm, mathrsfs, lscape}
\usepackage{verbatim}
\usepackage[all]{xy}
\usepackage{authblk}
\usepackage{color}

\usepackage{array, tabularx}

\usepackage{setspace}
\usepackage[top=1.5in, bottom=1.5in, left=1.3in, right=1.3in]{geometry}

\newtheorem{thm}{Theorem}[section]
\newtheorem{cor}[thm]{Corollary}
\newtheorem{lem}[thm]{Lemma}

\newtheorem{prop}[thm]{Proposition}
\newtheorem{prob}[thm]{Question}

\numberwithin{equation}{section}

\def \Q {\mathbb Q}

\def \C {\mathbb C}
\def \Z {\mathbb Z}
\def \P {\mathbb P}

\begin{document}
	
	\title{K\"ahler-Einstein metrics on $\Q$-smoothable Fano varieties, their moduli and some applications}
	\author{Cristiano Spotti\thanks{c.spotti@qgm.au.dk}}
	\affil{QGM, Aarhus University}

	\maketitle

	\begin{abstract} We survey recent results on the existence of K\"ahler-Einstein metrics on certain smoothable Fano varieties, focusing on the importance of such metrics in the construction of compact  algebraic moduli spaces of K-polystable Fano varieties. Moreover, we give some applications and we discuss some natural problems which deserve future investigations.    \end{abstract}
	
	

\section{Introduction}

Let $X$ be a smooth Fano manifold, i.e., a compact $n$-dimensional complex manifold with positive first Chern class or, equivalently, with ample anticanonical bundle $K_X^{-1}$.

In this survey we discuss the \emph{moduli problem} of this important class of complex varieties, showing its deep connections with the theory of \emph{canonical metrics} on complex manifolds. More precisely, we focus on the  so-called \emph{K\"ahler-Einstein} (KE) metrics. This correspondence can be thought as an higher dimensional generalization of the relations between the theory of compact complex curves and their natural algebraic degenerations to nodal curves (Deligne-Mumford moduli compactification), and  the theory of metrics with constant negative Gauss curvature and formation of hyperbolic cusps. However, crucially, in our higher dimensional KE Fano situation the value of the constant scalar curvature is \emph{positive}, fact that imposes, as we will see, important constraints on the possible degenerations of such spaces.

Beside a pioneering work of Mabuchi and Mukai in a special complex two dimensional case \cite{MM} (based on fundamental works on geometric limits of K\"ahler-Einstein manifolds in real dimension four by, among others, Anderson \cite{A} and Tian \cite{T90}), the precise picture on \textquotedblleft geometric compactified" moduli spaces for Fano manifolds remained unclear. In particular, it is important to note that the \textquotedblleft space" of all Fano manifolds is non-Hausdorff and Fano varieties may have continuous families of automorphisms. Thus some care has to be considered in studying such moduli problem.

However, the recent advances on the equivalence between existence of KE metrics and the purely algebro-geometric notion of K-stability on Fano manifolds \cite{CDS}, combined with the results \cite{DS14} by Donaldson and Sun on geometric limits of non-collapsing KE metrics (based on the so-called Cheeger-Colding-Tian theory of limit spaces), made possible to study in more detail \emph{moduli spaces of KE (or K-stable) Fano manifolds and their degenerations} showing that, if we restrict our attention to such special Fano varieties, their moduli theory becomes much well behaved. In particular, we now have a complete explicit picture in complex dimension two \cite{OSS}, and more general abstract results \cite{SSY, O15, LWX1, LWX2} in higher dimension regarding existence of weak KE metrics on singular Fanos and on the structure of the (compactified) moduli spaces.

Finally, for completeness, we should mention here in the introduction that the relation between canonical metrics and compact moduli spaces of varieties is also fundamental in the higher dimensional case  of varietes with negative first Chern class. Contrary to the Fano case,  we have that \emph{all} such smooth manifold admit KE metrics with negative  \textquotedblleft cosmological constant" by the works of Aubin \cite{Aub} and Yau \cite{Y}. If we consider singular varieties, K-stability is equivalent to   \textquotedblleft KSBA stability" \cite{O13, O11} (a condition on the singularities of canonical polarized varieties, generalizing Deligne-Mumford stability for curves, used to construct compact separeted moduli spaces,  e.g., \cite{Ale, K1, K2}).  Moreover,  it has been proved in \cite{BG} that  certain singular  KE metrics \cite{BG} always exist precisely on this type of singular varieties. We will briefly explain these relations in more detail at the end of Section $4$.

In conclusion, we can say that, at least for (anti)canonical polarized varieties, K-stability, with its relation with KE metrics, provides a \emph{unified} way to construct nice (compact) moduli spaces of algebraic varieties, and thus KE/K-moduli spaces are important objects to be further studied in the near future.

\subsection*{Acknowledgements}  This survey is an expanded version of a talk given at the INdAM meeting \textquotedblleft \emph{Complex and Symplectic Geometry}" held in Cortona, Arezzo (Italy), 12-18 June 2016. I would like to thank the organizers Daniele Angella, Paolo De Bartolomeis, Costantino Medori and Adriano Tomassini for the invitation, and Yuji  Odaka for comments on a draft of this note.  During the preparation of the  survey, the author has been partially supported by the AUFF Starting Grant 24285.

\section{An overview of Fano KE/K-moduli problem}

Let $X$ be an $n$-dimensional smooth Fano manifold and  let  $\chi(K_X^{-k})$ be the Euler characteristic of power of the anticanonical line bundle $K_X^{-1}$. By Kodaira's vanishing,  such Euler characteristic is equal to $h^0(K_X^{-k})$ and, moreover, it coincides with the Hilbert polynomial associated to the anticanonical polarization. 

Thus, for a given polynomial $h$, we can define the following \emph{moduli set}:
$$M_h:=\{X^n \mbox{ Fano  mfd   with} \, \chi(K_X^{-k})=h(k)\}/  \mbox{bi-holo}.$$
Being a \textquotedblleft parameter space" for certain algebraic manifolds, we would like this set to admit a  \emph{natural} algebraic/complex analytic structure of complex variety: i.e., if $\pi:\mathcal{X} \rightarrow B$ is a flat family where $\pi^{-1}(b)=X_b$ is a Fano manifold, the natural map $B \rightarrow M_h$ should be \emph{holomorphic} with respect to the analytic structure on $M_h$. 

However, for dimension $n\geq 3$ such structure cannot exist for trivial reasons, known as \textquotedblleft jumps of complex structures": there exist flat families of smooth Fanos $\pi:\mathcal{X} \rightarrow \Delta$  over the complex disc,  such that $X_t \cong X_s$ for any $t,s \neq 0$, but $X_0 \ncong X_t$, for $t \neq 0$. Thus $[X_0]\in \overline{[X_t]}$ (here the square bracket denotes the isomorphism class). Hence,  $[X_t]\in M_h$ would be a non-closed point, condition which is incompatible with the existence of a natural structure of complex analytic variety on $M_h$ (in particularly inducing a Hausdorff topology). A well-known concrete example of this phenomenon is given by deformations of Mukai-Umemura Fano 3-fold (see \cite{T97}, where the relations with KE metrics is discussed).

Thus the only hope to find a moduli space of Fanos which indeed admits a nice classical analytic structure is to \emph{restrict the class of Fanos to consider}. Of course this non-Hausdorff issue is typical in many moduli problems (e.g., moduli of vector bundles). An answer for solving this problem is usually found in restricting the attention to \textquotedblleft stable" vector bundles or, thanks to the Hitchin-Kobayashi correspondence, to bundles which admit Hermitian-Einstein metrics. 

This suggests that also \emph{in the case of varieties} we should look to certain \textquotedblleft stable" varieties or, somehow equivalently, varieties which carry special Riemannian metrics. However, understanding the \textquotedblleft right" stability condition to consider in the case of varieties turned out to be a highly non-trivial task, and many scholars worked in the last thirty years to better understand the relations between special metrics and algebraic stabilities, guided by the so-called \emph{Yau-Tian-Donaldson conjecture} (YTD for short): given a polarized complex manifold $(X,L)$, the existence of a K\"ahler metric with constant scalar curvature (cscK) in $2 \pi c_1(L)$ should be equivalent to certain purely algebraic notion of stability of $(X,L)$ (for a gentle introduction on this topic, focused on a moduli perspective, one can read \cite{RT}).  In particular, not all polarized manifolds carry canonical metrics, contrary to the Calabi-Yau or negative first Chern class case. Classical obstructions to the existence of such metrics are given by the reductivity of the automorphism group \cite{Mat} and the vanishing of the so-called Futaki invariant \cite{Fut}.

We are not going to describe the huge literature in the subject here, but instead we focus on our Fano case of anti-canonical polarized manifolds, where the natural differential geometric notion for a  canonical metric in $2 \pi c_1(K_X^{-1})$ reduces to the so-called \emph{K\"ahler-Einstein} (KE) condition. Recall that a KE metric on a Fano manifold is a K\"ahler metric $\omega \in 2 \pi c_1(K_X^{-1}) $ which satisfies the Einstein geometric PDE, necessarily with \textquotedblleft positive cosmological constant" (here normalized to $1$):
$$Ric(\omega)\left( =i\bar{\partial}\partial \log(\omega^n) \right)= \omega.$$ 
Thanks to the K\"ahler condition, such Einstein equation, in general obstructed, reduces to a complex Monge-Amp\`ere equation on a potential function and  thus it can be studied using techniques coming from pluri-potential theory.

In \cite{T97} Tian introduced the notion of \emph{K-stability}, extending the notion of Futaki invariant, stability condition later further generalized and made completely algebraic by Donaldson \cite{D02}. K-(poly)stability (be aware that sometimes people call K-stable something that for us is K-polystable!) is a \textquotedblleft   geometric invariant theory (GIT)-like"  notion of stability for varieties in which  one promotes an abstract version of the Hilbert-Mumford criterion as a definition for stability. 

A \emph{test-configuration} for $X$ (the analogous to a   one-parameter subgroup in standard GIT) is the datum of a $\C^{\ast}$-equivariant relative polarized normal flat family of schemes over $\C$: $$\C^\ast \curvearrowright  \left( (\mathcal{X},\mathcal{L}) \rightarrow \C \right),$$
such that over $1$ (hence away from zero), we have an isomorphism  $(X_1, \mathcal{L}_{|X_1})\cong (X, K_X^{-r})$. Clearly $\C^{\ast}$  acts naturally on the $d_k$-dimensional vector space $H^0(X_0,\mathcal{L}_{|X_0}^k)$ with weight $w_k$ on its top exterior power. Note that $X_0$ can be highly singular. 

The \emph{Donaldson-Futaki invariant} for the test configuration $(\mathcal{X},\mathcal{L})$ for $X$ (the analogous of weight in GIT) is the $k^{-1}$ coefficient in the Riemann-Roch expansion:
$$\frac{w_k}{kd_k}= C+ DF(X,(\mathcal{X},\mathcal{L}))k^{-1}+\mathcal{O}(k^{-2}).$$ $X$ is called \emph{K-stable} (rep. semistable) if and only if $DF(X,(\mathcal{X},\mathcal{L}))>0$ (resp. $\geq 0$) for all test-configurations, and \emph{K-polystable}  if and only if it is K-semistable and $DF=0$ iff $\mathcal{X}\cong X \times \C$.

Note that to define K-stability we haven't taken a specific embedding of $X$ in some \emph{fixed} projective space and  considered only embedded test configurations: the definition require a-priori to test stability for all possible equivariant degenerations inside any $\P^N$ where $X$ embeds, letting $N \rightarrow \infty$. For this reason testing K-stability from the actual definition is very challenging, even if we can reduce to the so-called \emph{special test configurations} \cite{LX}. However, some criteria related to the so-called  log-canonical-threshold, or to the very recent notion of Ding stability are available (e.g.,\cite{Fuj}).  More abstractly, K-stability may be thought as a GIT like notion on the stack of Fano varieties where the $DF$-invariant is actually realized as the weight of a stacky line bundle, called CM line bundle \cite{PT}. As we will see, this point of view is quite important for the moduli discussion.

We are now ready to state the fundamental theorem relating KE metrics with K-stability.

\begin{thm}\label{YTD}[YTD conjecture for Fano manifolds]. 
	Let $X$ be an $n$-dimensional smooth Fano manifold. Then 
	$$ \mbox{there exists a KE metric in } \, 2 \pi c_1(K_X^{-1})\,  \Longleftrightarrow \, X \, \mbox{ is K-polystable}.$$
\end{thm}

The direction \textquotedblleft$\Rightarrow$" has been proved in various degrees of generality by Tian \cite{T97}, Donaldson \cite{D08}, Stoppa \cite{St} and finally by Berman \cite{B}. The other direction is the content of the recent breakthrough of Chen, Donaldson and Sun \cite{CDS}. The proof uses a combination of analytic, geometric and algebraic techniques, in particular related to the notion of Gromov-Hausdorff (GH) convergence (notion that, as we will see, is deeply relevant also for the moduli problem). More recently different proofs have been found: via Aubin's continuity path \cite{DaS}, via K\"ahler-Ricci flow \cite{CSW}  or, for the case of finite automorphisms groups, via calculus  of variation  techniques \cite{BBJ}.

But now let us go back to the moduli discussion. We can define the \textquotedblleft differential geometric" KE moduli space of equivalence classes, up to \emph{biholomorphic isometries}, of KE Fano manifolds (with fixed Hilbert polynomial $h$):
$$ \mathcal{E} M_h:=\{(X,\omega)\, | \, \omega \, \mbox{KE} \}/ \sim.$$
Similarly, we can consider  the subset of $M_h$ defined by the algebro geometric condition of K-polystability:
 $$\mathcal{K}M_h:= \{[X] \, | \, X \, \mbox{K-ps}\} \subsetneq M_h.$$
Thus we have the following \emph{Hitchin-Kobayashi map for varieties}:
$$\phi_h: \mathcal{E} M_h \longrightarrow \mathcal{K} M_h,$$
naturally given by forgetting the metric structure, i.e.,  $\phi_h([(X,\omega)]):=[X]$. This map is:
\begin{itemize}
	\item well-defined, by \textquotedblleft$\Rightarrow$" in theorem \ref{YTD}.
	\item surjective, by \textquotedblleft$\Leftarrow$" in theorem \ref{YTD}.
	\item injective, by Bando-Mabuchi uniqueness \cite{BM}.
\end{itemize}

Thanks to the canonical metric structure  induced by the KE metric, we can now put a natural topology on the differential geometric moduli space $ \mathcal{E} M_h$. Such topology is essentially induced by the \emph{Gromov-Hausdorff (GH) distance} between compact metric spaces: given  two compact metric spaces, say $(S, d_S)$ and $(T, d_T)$, one defines
$$d_{GH}(S,T):= \inf_{S,T \hookrightarrow U} \inf  \{C >0 \, | \, S \subseteq N_C(T) \; \&  \; T \subseteq N_C(S) \},$$
where $N_C(S)$ denotes the distance $C$ neighborhood of $S$ isometrically embedded in a metric space $U$. The above defines a metric structure, in particular  a Hausdorff topology, on the space of isomorphisms classes of compact metric spaces. In practice, one usually estimates  the GH distance (which is sufficient for studying convergence) via maps $f:S \rightarrow T$ which are $\epsilon$-dense and $\epsilon$-isometries. See \cite{BB} for an introduction to such notion of convergence.

One of the immediate advantage of the GH topology is that, by its very definition, it gives a possible precise way to study \emph{degenerations} of Riemannian manifolds to \emph{singular} spaces.

Moreover, in our KE Fano case, we have to following remarkable pre-compactness property:  any sequence of complex  $n$-dimensional KE Fano manifolds $(X_i, \omega_i)$ subconverges  in the GH sense to a compact length metric space $S_{\infty}$ of  real Hausdorff dimension equal to $2n$. This follows by Gromov's theorem on convergence of Riemannian manifolds with Ricci uniformly bounded below and diameter bounded above (condition that in our case is implied by the positivity of the Ricci tensor, thanks to Myers'theorem) and by the volume non-collapsing condition, i.e., the volume of balls of radius $r$ is uniformly bounded below by $C r^{2n}$.

The metric space limit $S_\infty$ can be considered as a very weak limit. However, since we are considering limits of spaces which admit many additional structures (a Riemannian metric and a complex structure) it is natural to expect that in a suitable sense such structures are preserved in the limit. The first important result is the \textquotedblleft Riemannian regularity" provided by  
Cheeger-Colding theory \cite{CC} which shows that $S_\infty$ is actually an incomplete smooth Einstein space off a set of Hausdorff codimension $4$, and also  gives some geometric stratification of the singular set based on the local behavior of the metric structure (via splittings of metric tangent cones). The second regularity result is the recent theorem of Donaldson and Sun  \cite{DS14} which, in addition, shows that $S_\infty$  admits a natural structure of \emph{normal algebraic Fano variety}. Such structure is constructed by realizing the GH convergence as convergence of algebraic cycles in a sufficiently big projective space, via  \emph{uniform} Tian's $L^2$-orthonormal embeddings by plurianticanonical sections. Thus $S_\infty$ is homeomorphic to a Fano limit cycle $X_0$. Note that the above also gives a refinement of the GH topology, which now \textquotedblleft remembers" the complex structure (otherwise there is some ambiguity related to complex conjugations \cite{S}). More technically (see Section 4 for details), $X_0$ turned out to be a $\Q$-Gorenstein smoothable $\Q$-Fano variety, i.e., a \emph{normal} variety with $\Q$-Cartier anticanonical divisor and Kawamata-log-terminal (klt) singularities admitting nice smoothings and a weak KE metric.

It is important to remember that in complex dimension two the above convergence results were already  known by the works mentioned in the introduction of Anderson, Tian and others. In this situation, GH limits must have isolated orbifold singularities (which is precisely the klt condition in dimension two), i.e.,  quotients of $\C^2$ by finite subgroups of $U(2)$ acting freely on the $3$-sphere. Moreover, the KE metric is orbifold smooth, i.e., it extends to a smooth metric on the local orbifold covers.

Since, by the result of Berman \cite{B}, it is known that the direction KE implies K-polystability holds also for singular varieties, one can naturally define the  \emph{extended Hitchin-Kobayashi map}
$$\tilde{\phi}_h: \overline{\mathcal{E} M}^{GH}_h \longrightarrow  \overline{\mathcal{K}M}_h,$$
where $\overline{\mathcal{E} M}^{GH}_h$ is a compact Hausdorff topological space with respect to the refined GH topology obtained by adding all GH limits, and  $\overline{\mathcal{K}M}_h$ denotes the set of $\Q$-Gorenstein smoothable K-polystable $\Q$-Fano ($\Q$-smoothable for short) varieties up to isomorphism.

With all of this in mind,  it is natural to ask the following foundational questions:

\begin{enumerate}
	\item \emph{YTD for $\Q$-smoothable Fanos}: does any $X \in  \overline{\mathcal{K}M}_h$ admit a weak KE metric? (i.e., is $\tilde{\phi}_h$ surjective?)
	\item \emph{Existence of K-moduli}: does  $\overline{\mathcal{K}M}_h$ admit a natural algebraic structure such that $\tilde{\phi}_h$ is an homeomorphism with respect to the GH topology and the euclidean topology of the algebraic space?
	\item Can we find \emph{explicit examples} of such compactifications? 
\end{enumerate}

We start by discussing the dimension two case, i.e., the case of \emph{del Pezzo surfaces}. In this situation the answer to all such questions is complete. Note that in dimension one, if one does not consider weighted/cone angle case, the moduli problem for Fanos is clearly trivial, being $\P^1$ with the Fubini-Study metric the only such space.

\section{KE/K-moduli of del Pezzo surfaces}

In complex dimension two Fano manifolds are traditionally called \emph{del Pezzo surfaces}. Such varieties are also completely classified: they are given by $\P^2$, $\P^1\times \P^1$ and the blow-up of the plane in up to eight points in \textquotedblleft very general" position. Let us denote with $d=c_1^2(X)\leq 9$ their \emph{degree}  (which also uniquely determines the Hilbert polynomial, as a consequence of Riemann-Roch arguments). The problem of understanding which \emph{smooth} del Pezzo surface admits a KE metric was addressed in the seminal paper of Tian \cite{T90}. The answer is: they all admit such metrics, beside the well-known obstructed cases of the blow-up of the plane in one or two points. 

Since we are interested in moduli problem, we restrict to the case $d\leq 4$, i.e., to the case when there are non-trivial complex deformations. Together with recovering Tian's theorem in the smooth case and, somehow, providing a conceptual explanation why all smooth del Pezzo of degree $d\leq 4$ admit KE metrics, the following theorem computes via  explicit algebro-geometric techniques the GH compactifications of  such KE moduli spaces, classifying the geometric limits (which are KE del Pezzo orbifolds, by the result recalled in the previous section).

\begin{thm}[\cite{OSS}] \label{DP} For any positive  integer degree $d\leq 4$, there exists an \emph{explicit} compact algebraic space $\overline{M}_d^{ALG}$ (moduli space of certain degree $d$ del Pezzo orbifolds) such that the Hitchin-Kobayashi map
	$$\tilde{\phi}_d: \overline{\mathcal{E} M}^{GH}_d \longrightarrow   \overline{M}_d^{ALG} \left( \cong  \overline{\mathcal{K}M}_d \right)$$
	is a \emph{homeomorphism}, and $\mathcal{E} M_d$ is identified with a Zariski dense subset of $\overline{M}_d^{ALG}.$
\end{thm}

As we said in the introduction, the degree $d=4$ case was previously understood by Mabuchi and Mukai \cite{MM}.

An important differential geometric application, generalizing Tian's results in the smooth setting, is the following corollary, answering a conjecture of Cheltsov and Kosta \cite{CK}. Since KE del Pezzo orbifolds with orbifold groups contained in $SU(2)$ (i.e., with \emph{canonical singularities}) are classified and they always admit $\Q$-Gorenstein smoothings, from the above explicit KE/K-moduli compactification we have:

\begin{cor}[\cite{OSS}] KE del Pezzo orbifolds with orbifold groups at the singularities contained in  $SU(2)$ are classified.
\end{cor}

For example, KE del Pezzo orbifolds of degree three with such singuarities are precisely given by \emph{all} cubic surfaces in $\P^3$ with only nodal (i.e., $A_1$) singularities plus the toric cubic $\{xyz=t^3\}\cong \P^2/\Z_3$, since in this case the GH compactification agrees with the classical GIT quotient of cubic surfaces. Partial results where previously known (e.g.,  \cite{DT, CK, Shi}).
	
We now indicate the main passages in the proof of the above theorem  \ref{DP} for the interesting $d=2$ case.

Let $X_{\infty}$ be the GH limit of smooth degree $2$ del Pezzo surfaces (i.e., GH limit of KE double covers of $\P^2$ branched at a smooth quartics).
\begin{itemize}
	\item Step 1: we first improve our understanding of the singularity of $X_{\infty}$, combining Bishop-Gromov monotonicity formula (which shows that the order of the orbifold group at the singularity can be at most $6$) with the Koll\'ar and Shepherd-Barron classification of two dimensional $\Q$-Gorenstein smoothable quotient singularities \cite{KSB}.
	\item Step 2: using classification results of singular del Pezzo surface, we show that $X_\infty$ has to be given by the following hypersurfaces (of degree four and eight, respectively) in weighted projective spaces:
	\begin{enumerate}
		\item $X_\infty \cong \{f_4=t^2\}\subseteq \P(1,1,1,2)$;
	    \item $X_\infty \cong \{f_8=z^2+t^2\}\subseteq \P(1,1,4,4)$;
	\end{enumerate}
	\item Step 3: in both two cases there is a natural action of two groups on the parameter spaces: more precisely, in the first case $SL(3,\C)$ acts on the space of quartics $\P(Sym^4(\C^3))$ and in the second case $SL(2,\C)$ acts on $\P(Sym^8(\C^2))$. This gives two GIT quotients, with natural linearizations.
	\item Step 4: since $X_\infty$ is K-polystable by Berman \cite{B}, a comparison of stabilities using the CM line bundle, shows that $X_\infty$ has to be also stable with respect to the above classical notions of GIT stability. Next we can blow-up the first quotient semistable stack at the point corresponding to the double conic to get a space mapping to a categorical quotient $\overline{M}_2^{ALG}$: i.e., we can define $$\mathcal{M}_2:=[\P(Sym^4(\C^3))^{ss}/SL(3,\C)]\cup_{\{[q^2=t^2]\}}[\P(Sym^8(\C^2))^{ss}/SL(2,\C)]\rightarrow  \overline{M}_2^{ALG}.$$
	\item Step 5: since there exists at least one smooth degree two KE del Pezzo surface \cite{TY}, we can define a natural continuous map (the Hitchin-Kobayashi map) from $\overline{\mathcal{E} M}^{GH}_2$ to  $\overline{M}_2^{ALG}$. Finally a standard open-closed topological argument, combined with the fact that $\overline{M}_2^{ALG}$ is an Hausdorff space of del Pezzo orbifolds,  implies the statement.
\end{itemize}

Note that in general there are non-canonical singularities in the limits, e.g., the limit toric variety $X_\infty=\{x^4y^4=z^2+t^2\}\cong (\P^1\times \P^1)/\Z_4$ has in particular two singularities of type $\frac{1}{4}(1,1)$. This is related to the existence of \emph{torsion} Calabi-Yau ALE metric bubbles from limits of Einstein spaces \cite{Su}: loosing speaking, ALE bubbles are spaces that metrically model the formation of  singularities  in this non-collapsing setting, and thus they are somehow the equivalent of  hyperbolic collars in the (locally collapsing) curve case near the formation of a node.

 Let us also  note that, since an Einstein deformation of a smooth KE del Pezzo surface has to be KE \cite{LB},  such moduli spaces, quotenting with repect to the natural involution given by conjugating the complex structure, \emph{are explicit compactification of connected component of (real) Einstein moduli spaces} on the real oriented manifolds $\C \P^2  \sharp k \overline{\C \P^2}$ with $5\leq k\leq 8$. See, for example, the general discussion on Einsten moduli by Koiso \cite{Koi}. The topological types of our spaces could be easily understood.

 In degree one the construction of the algebraic compactification is more involved, since bi-meromorphic contractions are used (which also cause the a-priori loss of projectivity for the moduli space). However, in all cases the algebraic compactifications obtained show that such moduli spaces admit \textquotedblleft more structure": namely, they are \emph{KE moduli Artin stack} (essentially \'etale covered by affine GIT quotients parametrizing del Pezzo orbifold deformations, see definitions $3.13-14$ in \cite{OSS}). This is related to the Alper's notion of \emph{Good Moduli Spaces} for an Artin stack \cite{Al}.
 
 Finally, the fact that our explicit moduli agrees with the K-compactification (that is, all $\Q$-smoothable K-polystable del Pezzo surfaces appear in their boundary) is a consequence of more general results which we are going to discuss in the next section.

\section{KE metrics on $\Q$-smoothable Fano varieties}

By what we explained in Section $2$, in order to study the boundary of the moduli problem is natural to consider  $\Q$-Gorenstein smoothings of a K-polystable $\Q$-Fano ($\Q$-smoothings for short) varieties, i.e., flat families $\mathcal{X}\rightarrow \Delta$ over the complex disc where:
\begin{itemize}
	\item $X_0$ is a normal (K-polystable) Fano variety with $\Q$-Cartier canonical divisor (i.e., some power is a line bundle) satisfying $K_{\hat{X}_0/X_0}=_{\Q}\sum_i a_i E_i$, with $a_i>-1$, for any log-resolution $\hat{X}_0$.  Equivalently, from a more differential geometric view-point, for any $p\in X_0$,  $\int_{U\cap X_0^{reg}} s^{\frac{1}{m}} \wedge \bar{s}^{\frac{1}{m}} < \infty$, where $s$ is a local trivialization over a small neighborhood $U$ of the $m$-th power of the canonical bundle $K_{X_0}^m$ near $p$.
	\item $K_{\mathcal{X}/ \Delta}$ is $\Q$-Cartier. 
	\item $X_t$ is smooth.
\end{itemize}
We remark that \emph{not} $\Q$-\emph{Gorenstein} smoothings of $\Q$-Fano varieties exist, but such deformations are not relevant for KE/K-moduli problems. We can now state the main theorem of this section, which in particular give an answer to question one in Section 2.
\begin{thm}[\cite{SSY}]\label{QF}
	Let $\mathcal{X}\rightarrow \Delta$ be a $\Q$-smoothing of a K-polystable variety  $X_0$. Then $X_0$ admits a weak KE metric $\omega_0$. Moreover, if $Aut(X_0)$ is finite, $X_t$ admit smooth KE metrics $\omega_t$ for $t$ sufficiently small, and $(X_t,\omega_t) \rightarrow (X_0, \omega_0)$ in the GH topology.
\end{thm}
Thanks to the above correspondence between metric limits and flat families, as a corollary we have the following algebraic separatedness statement.
\begin{cor}[\cite{SSY}]
	If two $\Q$-smoothings of K-stable $\Q$-Fanos (with finite automorphisms) agree away from the singular fiber, then the central (singular) fibers are isomorphic.
\end{cor}

Before giving a quick survey of the main ideas in the proof, it is useful to discuss some properties of weak KE metrics, which are the natural pluri-potential theoretic generalization on singular varieties of smooth KE metrics (\cite{EGZ}, or \cite{Dem} for a recent survey).

Near any point $p \in X_0$ such weak KE metrics are given by the restriction of the $i \partial \bar \partial$ of a \emph{continuous} potential for an embedding in $\C^N$ of the analytic germ of the singularity. As two dimensional orbifold singularities show, such regularity  for the potential is essentially optimal. Regarding more geometric considerations, we have that the regular part $X_0^{reg}$ is a \emph{smooth} incomplete KE space, and its metric completion 
$\overline{X_0^{reg}}$ topologically agrees with $X_0$. 

The actual \textquotedblleft asymptotic behavior" at the singularities of these weak KE metrics is quite delicate, and  more complicated with respect to the two dimensional orbifold case. Recently we have seen important results which put some light on the aspect of the metric near the singular locus,   at least when we consider singular KE spaces arising as limits of smooth ones (note that  orbifold singularities appear only exactly in complex codimension two, by the famous Schlessinger's rigidity of quotient singularities). From a metric measure theoretic perspective,  it is known that 
the metric \textquotedblleft looks the same" at all sufficiently small scales near a singularity (uniqueness of metric Calabi-Yau tangent cone \cite{DS15}). But, as first observed by Hein and Naber \cite{HN}, phenomena of local jumping of complex structures can happen when \textquotedblleft zooming" to find such metric tangent cone. For example, it is expected that in complex dimension $3$, metric tangent cones at the isolated singularities of type $A_k$ (i.e., locally analytical of type $x_1^2+x_2^2+x_3^2+x_4^{k+1}=0$) for $k \geq 3$ should all be isometric to the flat cone $\C \times \C^2/\Z_2$ (singular along a line). See \cite{D16} for a discussion.

Roughly speaking (but the situation is slightly more subtle in reality), these jumping phenomena have their origins in the fact that  typical complex links of klt singularities are  Fano varieties, but the existence of a Calabi-Yau cone metric model implies that such links have to be KE! So there is indication that some  notion of \emph{stability for singularities} is required (related to Sasaki-Einstein stability, in the simplest cases). Some recent works, such as \cite{Li16, LL}, are trying to understand this picture from an algebraic perspective. However, in certain situations (e.g., for the  $A_1$ case, where a CY cone metric can be found via Calabi's ansatz \cite{C} and where the corresponding smoothing bubble was explicitly found by Stenzel \cite{Ste}) it is expected, and very recently proved in the CY case by Hein and Sun \cite{HS}, that the weak KE metrics are \emph{polynomially asymptotic} to the CY cone models in a suitable local holomorphic gauge.

Let us now briefly describe the strategy in the prove of theorem \ref{QF}. The very rough idea consists in constructing the weak KE metric on the singular fiber via a GH limit of certain conically singular KE metrics on the smooth nearby fibers, thus running in \emph{families} the so-called \textquotedblleft Donaldson's cone angle path", i.e., the continuity path used by \cite{CDS} for proving the YTD conjecture  in the Fano case.

\begin{itemize}
	\item Step 1: by a Bertini's type argument we can take a divisor in $\mathcal{D} \in |-\lambda K_{\mathcal{X}/\Delta}|$, for $\lambda$ big enough,  which gives a smooth pair $(X_t,D_t)$ when restricted at $t\neq 0$ and a klt pair $(X_0,(1-\beta)D_0)$ on the singular fiber for $\beta \leq 1$. Thus we want to consider the following  \emph{two parameters family} of PDEs:
	$$Ric(\omega_t+i\partial\bar{\partial} \phi_{t,\beta})=(1-(1-\beta)\lambda) (\omega_t+i\partial\bar{\partial} \phi_{t,\beta})+2\pi (1-\beta)\delta_{D_t},$$
	where  $\omega_t$ is the restriction of a Fubini-Study metric from an embedding of the family and $\delta_{D_t}$ the current of integration along $D_t$.  At least when $t\neq0$ a solution of this equation is a KE metric with cone angle equal to $\beta$  along $D_t$ (e.g. \cite{JMR} or \cite{GP}).
	\item Step 2: via a log-canonical-threshold argument one shows that the above equation has a positive KE solution for all $t$, if $\beta$ is sufficiently small.
	\item Step 3: using some pluri-potential techniques (e.g., Berndsson's positivity of direct images) one can find, for fixed $\beta$, an a-priori bound of type $||\varphi_{t,\beta}||_{L^{\infty}}\leq C(\beta),$ for $t\neq 0$. Taking the universal embedding in $\P^N$ provided by the conical generalization of Donaldson and Sun convergence theorem \cite{CDS} III, one sees that the conical metrics on the smooth fibers GH converge to the weak conical KE metric on the central fiber.
	\item Step 4: the above convergence is used to prove that the function 
	$$\beta_t:=\sup \{ \beta \in (1-\lambda^{-1}, 1] \, |\, \exists \, \omega_{t,\beta} \, \mbox{KE on } X_t \},$$
	is a lower semi-continuous function in the euclidean topology of the disc. This is connected with some properties of the automorphism groups.
	\item Step 5: the above semi-continuity result, combined with a gap argument for some natural energy functional (Aubin's energy), gives that the set of cone angles $\beta$s for which a weak conical KE metric exists on $X_0$ is open. The closeness follows again by taking limits from the smooth nearby fibers.
	
\end{itemize}

Thus the KE metric on $X_0$ is constructed as a kind of  \textquotedblleft diagonal" GH limit of cone angle KE metrics $\omega_{t,\beta(t)}$ with $\beta(t)\rightarrow 1$ as $t\rightarrow 0$. In particular, it is a weak KE metric thanks to the regularity theory for GH limits.

\subsection{Algebraic structure on Fano KE/K-moduli}

Theorem \ref{QF} above shows that the YTD conjecure also holds in the case of $\Q$-smoothable Fano varieties and it provides a natural correspondence between flat limits and GH convergence, at least in the case of finite automorphisms groups. The next step is then related to the construction of a natural algebraic structure on the differential geometric KE moduli space or, equivalently, showing that our question two asked at the end of the second section admits a positive answer. 

The rough idea for constructing such algebraic moduli space of $\Q$-smoothable K-polystable Fano varieties is the following. One knows that, being KE, such varieties have linear reductive automorphism groups. Thus one can consider a Luna's slice type argument applied to the Hilbert scheme in the uniform $\P^N$ embedding where all GH limits of smooth KE spaces (say with  fixed Hilbert polynomial $h$) live. Here one shows that, \'etale locally, K-polystability is completely captured by some \textquotedblleft local GIT" stability on a small enough (affine) slice. The actual argument is similar to Step 4-5 in the construction of explicit moduli space of del Pezzo surfaces. Thus this expected  \textquotedblleft local GIT picture" (e.g. \cite{S}), generalizing the one in the smooth case obtained by Broennle  \cite{Br} and Sz\'ekelyhidi \cite{Sz}, provides the natural, compatible with the GH topology, \emph{algebraic atlas for a compact moduli space}. The above has been fully proved by Li, Wang, Xu \cite{LWX1} and, independently, by Odaka \cite{O15} (using theorem   \ref{QF} recalled above and some preliminary propositions in the first version of \cite{LWX1}), generalizing \cite{O132}. In conclusion we have: 

\begin{thm}[\cite{LWX1},\cite{O15}]
In any dimension, $\overline{\mathcal{K}M}_h$ admits a natural algebraic structure (given \'etale locally by affine GIT quotients) such that the map $\tilde{\phi}_h$ is a homeomorphism.
\end{thm}

As in the del Pezzo case, this moduli space carries more structure: in particular, $\overline{\mathcal{K}M}_h$ is a categorical quotient of a  KE/K-moduli stacks as we have previously discussed. Moreover, further analyzing the structure of such moduli spaces (in particular studying openness of K-semistability), the authors in \cite{LWX1} showed that 
$\overline{\mathcal{K} M}_h$ is \textquotedblleft dominated" by a good Artin moduli stack $\overline{\mathcal{K M}}_h$ of \emph{K-semistable} $\Q$-smoothable Fano varieties, with a unique K-polystable point in the K-semistable equivalence classes. However, as we have seen, we stress that at present the existence/construction of such algebraic moduli spaces of $\Q$-smoothable K-polystable Fano varieties depends crucially on  transcendental complex analytic techniques related to KE metrics. Some discussion on the potential dependence on $N$ (the dimension of the projective space were all GH limits live) of the algebraic structure on the compactified moduli space can be found in the original papers. Moreover, it will be very important to find a purely algebraic way to form such moduli spaces and, furthermore, to remove the smoothability hypothesis used in the present construction. We expect new birational geometric techniques to be relevant for this progress.

Finally, we mention that such KE/K compact moduli spaces of Fano varieties are the analogous of the KSBA compactification of moduli spaces of manifolds with negative first Chern class \cite{K2}, and thus a special instance of the more general theme: relations between special K\"ahler metrics and moduli of polarized varieties. In \cite{O13, O11} Odaka showed that K-stability is equivalent to the KSBA conditions on the singularities (semi-log-canonical) of a variety (satisfying the conditions G$1$ and S$2$) with ample canonical divisor  required to form compact moduli spaces. Moreover, Berman and Guenancia showed in \cite{BG} that precisely on varieties with such type of singularities is possible to construct weak KE metrics of negative scalar curvature. Here the metric can be \emph{complete} near the non-klt locus, and locally collapsing (this is precisely the higher dimensional analogous of the hyperbolic cusps in the complement of a node in a DM stable curve). Thus, even if indirectly, one recovers the equivalence between K-stability and (negative) KE metric. However, the complete metric  convergence picture is not fully understood, due to these collapsing phenomena. It is known that if $X_0$, the central fiber of a smoothing, has only simple normal crossing singularities, then the KE metrics in the nearby fiber naturally converge to complete KE metrics on the irreducible components of $X_0$ (e.g., \cite{T93}, \cite{R}). Some properties of (special) collapsing regions has been recently studied by Zhang in \cite{Z1} inspired by the SYZ picture in collapsings of Calabi-Yau manifolds. Related to the last point, we should mention that canonical metrics should be relevant also (at least) in the study of compactified moduli space of polarized Calabi-Yau manifolds: the non-collapsing case is well understood (e.g., \cite{Z2}), however the full GH collapsing to lower dimensional spaces (e.g., \cite{GTZ}) remains quite mysterious (but see conjectures of  Gross and Wilson, and Kontsevich and Soibelman, e.g., \cite{KS}, where collapsing to certain spaces of real dimension at most equal to half of the original dimension is expected), and possibly related to certain moduli of tropical varieties \cite{O14}. For relations with algebraic geometry, fixing the polarization in the study of degenerations of Calabi-Yau manifolds is going to be essential, as the purely trascendental collapsing of K3 surfaces to real three dimensional spaces in \cite{Fos} suggests.  

\section{Some applications and future perspectives}

In this last section we describe some possible applications of the previously discussed results and, moreover, we will mention some natural problems to be considered in the near future.

The first application, more differential geometric in nature, consists in using singular KE metrics to construct examples of \emph{smooth K\"ahler metrics of constant scalar curvature} (cscK), a notorious difficult problem, via certain geometric transitions. Next we discuss some properties related to the study of the \textquotedblleft geometry" of KE/K-moduli spaces or stacks. Finally, we briefly mention the problem of understanding \emph{explict} examples of KE/K-moduli spaces.

\subsection{Generalized cscK conifold transitions}

Through this section let $X_0 \hookrightarrow \mathcal{X}\rightarrow \Delta$ be a $\Q$-smoothing of a K-stable Fano variety (with discrete automorphism group). By theorem  \ref{QF}, $X_0$ and its sufficiently small deformations $X_t$ are KE, and moreover the family is continuous in the GH topology. Now let us take a \emph{resolution}  $\hat{X}_0$ of the singular variety $X_0$.

\begin{prob}
	Can we find a family $g_\epsilon$ of \textquotedblleft canonical" K\"ahler metrics on a resolution $\hat{X}_0$ which also degenerate, as $\epsilon \rightarrow 0$, to the singular KE space $X_0$?
\end{prob}

The natural notion of best metric to consider on the resolution $\hat{X}_0$ is given by the more general notion of cscK metrics. In a loose sense, we can think of KE metrics on the smoothings to be a family of metrics where the underlying symplectic structure is fixed while the complex structures changes and becomes degenerate.  For the metrics on the resolution the relations is the opposite one: the complex structure is now fixed but the symplectic structures vary (the $\epsilon$ parameter being related to let the K\"ahler classes of the metrics approaching a special point in the boundary of the K\"ahler cone on $\hat{X}_0$). We call such paths of canonical metrics connecting in the GH sense smooth complex manifolds, in general not diffeomorphic, through a singular variety   \emph{generalized cscK conifold transitions}. Such terminology originates from a similar geometric situation considered in  Physics for Calabi-Yau 3-folds. 

The construction of these geometric transitions is expected to be hard in general. However, in the case when $X_0$ has only isolated singularities of some special type, one can hope to show existence of cscK metrics on some resolutions via gluing techniques, similar to the strategy used in \cite{AP} in the case of (orbifold) smooth metrics. For example, if the singularities of $X_0$ are locally analytically modeled on the blow down of the zero section of the canonical bundle of a KE Fano manifold (in general not orbifold) and the KE metric on  $X_0$ is asymptotic near the singularities to the conical CY cone metric given by the Calabi's ansatz \cite{C}, we can prove the following:

\begin{thm}[\cite{AS}]
	Under the above hypothesis, $X_0$ has a natural crepant  resolution $\hat{X}_0$  admitting a family of cscK metrics of positive scalar curvature converging to the KE metric on $X_0$ in the GH topology, and thus $X_0$ is the degenerate variety of a generalized cscK conifold transition.
\end{thm}

The above theorem is a special case of  more general results in  \cite{AS} (combined with theorem \ref{QF}), where $X_0$ is not assumed to be Fano (e.g., it could have a KE metric of zero or negative Einstein constant), nor smoothable, and the singularities belong to a bigger class. We crucially remark that having the needed asymptotic behavior of the weak KE metric near the singularities is in general a major problem. However, as we have previously recalled, for the smoothable Ricci-flat case (but modifications of the arguments should also work for KE metrics with different sign of the Einstein constant) Hein and Sun have recently shown \cite{HS} that the required  asymptotic decay property of the weak KE  metric for isolated singularities of the type considered in the above theorem. 

\subsection{Geometry of  KE/K-moduli stacks}

As we have explained in the previous sections, inside the non-separated, non-proper moduli stack $\mathcal{M}_h$ of Fano varieties, in general not of finite type even if we restrict the attention to $\Q$-smoothable Fanos, we can find a nice subspace of K-semistable $\Q$-smoothable Fano varieties $\overline{\mathcal{K M}}_h$ mapping canonically to a compact algebraic space $\overline{\mathcal{K}M}_h$, the coarse moduli variety of smoothable K-polystable or KE objects. Thus we ended up in a set-up of good moduli space for an Artin stack \cite{Al}, here described by a kind of generalized GIT quotient space with respect to a more abstract stability notion (K-stability) given by the CM line bundle on the moduli stacks $\mathcal{M}_h$  (see \cite{H-L} and  \cite{H} for more information on this view point, which we think will be relevant in the future).
The natural next step in this moduli theory is to understand further the geometry of such  spaces. For example, even if the CM line bundle  $\lambda_{CM}$ is in general not ample \cite{FR}, it is expected -thanks to its relation with Weil-Petersson geometry \cite{FS}- that, once descended to the coarse algebraic space, it becomes a natural $\Q$-polarization for $\overline{\mathcal{K}M}_h$, which then will be a \emph{projective variety} (the quasi-projectivity of the part parameterizing smooth KE/K-polystable Fanos has been recently shown in \cite{LWX2}).

Furthermore, it would be interesting to study properties of \textquotedblleft subvarieties" of the moduli stacks (i.e., families of K-stable Fano varieties), canonical bundles or sheaves on such moduli spaces, and cohomological properties of them, similar to the ones studied for curves and varieties with ample canonical class.  

As a toy example of such possible investigations, we will now compute the CM-volume of a simple curve (i.e., a family over a one dimensional space) in the KE/K-moduli space of degree $3$ del Pezzos, i.e., cubic surfaces. We first need an easy lemma:
\begin{lem} Let $\gamma: \mathcal{X} \rightarrow \mathcal{C}_g$ be a  curve of degree $d$ del Pezzo orbifolds with generically smooth fibers for which $K^{-1}_{\mathcal{X}/\mathcal{C}}$ makes sense. Then
	$$c_1(\lambda_{CM}(\mathcal{X} \rightarrow \mathcal{C}_g) )=6d(1-g) -c_1^{3}(\mathcal{X}).$$
\end{lem}
 In fact, by the definition of the CM line bundle for the relative anticanonical polarization and by Grothendieck-Riemann-Roch, we have $c_1(\lambda_{CM}(\mathcal{C}))=-\gamma_\ast ( c_1^3(K^{-1}_{\mathcal{X}/\mathcal{C}}))$. Hence, $c_1(\lambda_{CM}(\mathcal{C}))= -\gamma_\ast \left( c_1^3(K^{-1}_{\mathcal{X}})+3c_1^2(K^{-1}_{\mathcal{X}})c_1(\gamma^\ast K_{\mathcal{C}})\right)$ which is indeed equal to $-c_1^{3}(\mathcal{X})-6d(g-1)$, as claimed.

Thus, for example, if $d=3$ and $\mathcal{C}=\P^1$, we have $c_1(\lambda_{CM})=18-c_1^3(\mathcal{X}).$ Moreover note that, by the positivity of the CM line bundle, if $\mathcal{X} \rightarrow \mathcal{C}_g$ is a  \textquotedblleft K-polystable curve" then $c_1^3(\mathcal{X})\leq 6d(1-g)$. Similar Chern numbers inequalities can be founded in higher dimension too.

Now, if we take a \emph{generic pencil} of cubic surfaces $tc_1+sc_2=0$,  by genericity we may assume  that the generic member in the associated \emph{Lefschetz's fibration} $\mathcal{X}\rightarrow \P^1$  is smooth and the singular fibers have only one nodal $A_1$-singularity.  Thus we have the following \textquotedblleft intersection number computation".
\begin{prop}The degree of the CM line bundle on the base of a generic Lefschetz's fibration  of (K-stable by theorem \ref{DP}) cubic surfaces is equal to $$c_1(\lambda_{CM}(\mathcal{X}\rightarrow \P^1))=8.$$
\end{prop}
	For this, thanks to the previous lemma,  it is sufficient to compute $c_1^3(\mathcal{X})$, where $\mathcal{X}=Bl_{\Sigma_g}\P^3$ with $\Sigma_g=c_1\cap c_2$ surface of genus $g=10$, by adjunction. Since $-K_{\mathcal{X}}=4H-E$, where $H$ is the pull-back of the hyperplane bundle of $\P^3$ and $E$ is the exceptional divisor ($\P^1$ bundle over $\Sigma_g$), we have that
	$c_1^3(\mathcal{X})=64H^3-48H^2E+12HE^2-E^3.$
	But $H^2E=0$, $HE^2=9E.f=-9$ (where $f$ is a fiber of the $\P^1$-bundle), and $-E^3=N_{\Sigma_g}=-K_{\P^3}.C+2g-2=54$. Hence $c_1^3(\mathcal{X})=10$, which implies the result. 
	
We expect  similar computations to be  relevant in the study of  properties of the Picard group of  K-moduli stacks $\overline{\mathcal{K M}}_h$.

\subsection{Examples of Fano KE/K-moduli}

Beside the complex dimensional two case, (where  the KE/K-moduli picture is complete, at least for the components corresponding to compactifications of smooth surfaces), in higher dimension we are completely lacking of explicit examples. There are several reasons to look for such examples. For us the two more important ones are:
\begin{itemize}
	\item they will provide a complete understanding of which Fano manifolds in a given family admit KE metrics.
	\item they may provide hints to study recurrent properties of K/KE-moduli spaces.
\end{itemize}
We expect that the techniques developed in the proofs of the theorems presented and discussed in this survey note (e.g., stability comparisons, local moduli picture, properties of singularities, etc...) will be essential in the future studies. Natural situations to investigate are given by Fano $3$-folds, log settings, \textquotedblleft special" Fanos, non-smoothable KE del Pezzo orbifolds. It is natural to believe that explicit K-moduli compactifications could be found by birational modifications of standard GIT quotients, as we have shown for the two dimensional del Pezzo case.

\end{document}